\documentclass{article}
\usepackage{graphicx} 
\usepackage{amsthm}
\usepackage{amsmath}
\usepackage{amssymb}
\usepackage[utf8]{inputenc}
\usepackage[english]{babel}
\usepackage{amsthm}
\usepackage{biblatex}
\renewbibmacro{in:}{}
\newtheorem{theorem}{Theorem}[section]
\newtheorem{lemma}[theorem]{Lemma}

\theoremstyle{definition}
\newtheorem{definition}{Definition}[section]

\theoremstyle{remark}

\title{REU - Schwartz Functions and Compactifications}
\author{Jonah Marcus, Molly Sager \\ under supervision of Ahmad Reza Haj Saeedi Sadegh }
\date{May 2023}

\begin{document}

\maketitle

\begin{abstract}
    In the early 20th century, Laurent Schwartz observed that we can identify functions that extend smoothly to the point at infinity of one-point compactifications of Euclidean spaces. We show a similar result for a different compactification of Euclidean spaces, namely, the real projective spaces.
\end{abstract}

\section{Introduction}
The one-point compactification of topological spaces has been an established concept since the early 20th century. Pavel Alexandroff showed that one can adjoin a single point to a topological space such that the resulting space is compact. Later in the century, the French mathematician Laurent Schwartz showed that smooth functions with rapid decay at infinity, now known as Schwartz functions, can be identified as those functions that extend smoothly to the infinity point of the one-point compactification of $\mathbb{R}^n$. Our aim is to show that Schwartz functions extend smoothly to a flat function at the infinity points of the real projective space, $\mathbb{RP}^n$, considered as a compactification of $\mathbb{R}^n$. \\

\quad We begin this paper by introducing all preliminary concepts necessary for obtaining our result. We proceed to the construction of the real projective space and justify its manifold structure. Then, we define Schwartz functions and the notion of smooth manifolds. Following this, we will show Schwartz' result for the one-point compactification of $\mathbb{R}^n$, and extend it to a generalization of this notion for real projective spaces.

\section{Preliminaries} 
In this section we introduce all the necessary concepts needed to understand real projective spaces as an abstract topological manifold. 

\subsection{Topological Spaces}
We begin with a discussion of \emph{topological spaces} and some of the important notions concerning topological spaces: 

\begin{definition}[Topological Space]
    A \emph{topological space} on a set X is a pair $(X,\tau)$ where X is the set itself and $\tau$ is a collection of "open" subsets of X satisfying the following conditions:
    \begin{enumerate}
        \item $\emptyset , X \in \tau$ 
        \item $C \subset \tau \Rightarrow \bigsqcup_{U \in C}U \in \tau$
        \item $C_1,...,C_n \in \tau \Rightarrow \bigcap_{i=1}^{n}C_i \in \tau$
    \end{enumerate}   
\end{definition}

A number of examples of topological spaces come to mind. First, consider the set $\mathbb{R}^n$. We may define a topology on $\mathbb{R}^n$ by letting $\tau$ be the class of all sets generated by the arbitrary union and finite intersection of all open balls in $\mathbb{R}^n$. We call this the \emph{standard topology} of $\mathbb{R}^n$. Another example of a topological space would be $S^n$. To see that this is the case one may notice that $S^{n}\subset \mathbb{R}^{n+1}$. With this observation we may define the topology on $S^n$ as $\{x \cap S^n : x\in \tau_R \}$, where $\tau_{R}$ is the standard topology of $\mathbb{R}^{n+1}$. We may generalize this observation in the following way. Suppose $(X,\tau)$ is a topological space and $Y\subseteq X$. We may define a topology on Y by letting the topological space be defined as $(Y, \{x \cap Y : x\in \tau \})$. We call this the \emph{subset topology}. 

\begin{definition}[Second-countability]
    A topological space $(X,\tau)$ is said to be \emph{second-countable} if every $x \in \tau$ may be generated from a countable subset $\mathcal{B} \subset \tau$ through arbitrary unions of elements of $\mathcal{B}$
\end{definition}

\begin{definition}[Hausdorff criterion]
    A topological space $(X,\tau)$ is said to be \emph{Hausdorff} if for every $x, y \in X$, there exists open subsets $U, V \in \tau$ such that $U\cup V = \emptyset$ with $x \in U$ and $y \in V$.
\end{definition}

One can see that the metric topology on $\mathbb{R}^n$ ensures that this space is Hausdorff and second-countable. One can see this by letting the topology of this space be generated by the open balls in $\mathbb{R}^n$ that are centered at rational coordinates with rational radii. This set of open balls will be denoted as $B_{\mathbb{Q}}^n=\{ B_r(x):r\in \mathbb{Q}, \: x \in \mathbb{Q}^n \}$, where $B_r(x)$ denotes the open ball centered at x with radius r. The topological space generated from these balls is then defined as

\begin{center}
   $\mathbb{R}^n_{\mathbb{Q}}=(\mathbb{R}^n,\{\bigsqcup_{U \in C}U:C \subset B_{\mathbb{Q}}^n\})$
\end{center}

 To show this space is Hausdorff we consider two points x and y in $\mathbb{R}^n$ where $x \neq y$. It is clear that $B_{\frac{|x-y|}{2}}(x)\cap B_{\frac{|x-y|}{2}}(y)=\emptyset$ and that $B_{\frac{|x-y|}{2}}(x),B_{\frac{|x-y|}{2}}(y) \in B_{\mathbb{Q}}^n$. From this we may conclude that $\mathbb{R}^n_{\mathbb{Q}}$ is Hausdorff. The second-countability of this space is ensured by the fact that it is generated by the elements of $B_{\mathbb{Q}}^n$. Each element of this set is uniquely determined by an (2n+2)-tuple of integers. This is because the radius is determined by 2 integers and each element of $\mathbb{Q}^n$ is determined by 2 integers for each of the n coordinates. This set is countable.

\begin{definition}[Compactness]
    A collection of open sets C is said to be a \emph{open cover} of a set X if $\bigcup_{x \in C}x = X$. A collection of sets D is said to be a \emph{open subcover} of C if $D\subseteq C$ and $\bigcup_{x \in D}x = X$. A topological space, say T, is said to be \emph{compact} if for every cover of T there also exists a finite open subcover.
\end{definition}

$S^n$ is an example of a compact space. To see this we will use the Heine-Borel theorem. This theorem says that a subset of euclidean space is compact if and only if it is both closed and bounded. To show that $S^n$ is bounded we use the euclidean norm $|\:|_n:\mathbb{R}^n\rightarrow [0,\infty)$. One can see that the image of $S^n$ under this norm is $|\:|_{n+1}(S^n)=\{1\}$. This means $S^n$ is bounded. We know $S^n$ is closed Since $\overline{S^n}=\mathbb{R}^n \setminus S^n$ is open. Since $S^n$ is bounded and closed the Heine-Borel theorem tells us that this space is also compact. But one may also notice that since $S^n \subset \mathbb{R}^{n+1}$ that $S^n$ has the subset topology. This implies that $S^n$ is not only compact but also Hausdorff and second-countable as it inherits these properties from $\mathbb{R}_{\mathbb{Q}}^{n+1}$.

\begin{definition}[Continuity]
    We say a map between two topological spaces $f:(X,\tau_1) \rightarrow (Y,\tau_2)$ is \emph{continuous} if for every $C \in \tau_2$ we have $\tau_1$, $f^{-1}(C) \in \tau_1$.
\end{definition}

This notion of continuity will be used in our discussion of smooth abstract manifolds so a discussion of it will be omitted until then. 

\subsection{Manifolds}
With this discussion of topological spaces concluded we may introduce the concept of homeomorphisms and diffeomorphisms. This will allow us to begin an exploration of manifolds. 

\begin{definition}[Smoothness]
    A function f is \emph{smooth} on some open subset of $\mathbb{R}^n$ if every partial derivative exists at every point in said domain to infinite order. We can say a function f is smooth on a subset of $\mathbb{R}^n$ if there exists an extension of f to an open subset and f is smooth on this open subset. 
\end{definition}

There are many examples of smooth functions, including polynomials, exponential functions, analytic functions, and many others. In fact, one needs to take effort to produce a non-smooth function.

\begin{definition}[Homeomorphism]
    Consider two topological spaces $(X,\tau_1)$ and $(Y,\tau_2)$. A \emph{homeomorphism} is a bijective and continuous map $f:X \rightarrow Y$ where $f^{-1}$ is also continuous. 
\end{definition}

\begin{definition}[Diffeomorphism] 
    A function $f: M \to N$ where $M,N$ are manifolds is called a \emph{diffeomorphism} if $f$ is a smooth bijection and $f^{-1}$ is also smooth.
\end{definition}

In the definition of a diffeomorphism we used the concept of a \emph{manifold}. Intuitively a manifold can be thought of as a space that is locally euclidean. This means a diffeomorphism can be understood as a continuous "stretching" of one manifold into another. To truly understand these concepts we must have a formal definition of manifolds. We first consider manifolds explicitly embedded in $\mathbb{R}^n$ and then focus on abstract topological manifolds. 

\begin{definition}[Smooth Real Manifolds]
    A subset $M\subseteq \mathbb{R}^n$ is said to be a \emph{real smooth manifold of dimension m} if each $x \in M$ has a neighborhood that is diffeomorphic to to an open subset of $\mathbb{R}^m$. 
\end{definition}

$S^n$ is an example of a smooth n-dimensional real manifold. To show this is the case we use half spheres as covers of $S^n$ and show these half spheres are diffeomorphic to an open subset of $\mathbb{R}^n$. We know $S^n=\{(x_1,x_2,...,x_{n+1}):x_1^2+x_2^2+...+x_{n+1}^2=1\}$. To use half spheres as covers of $S^n$ we define them as:

\begin{center}
    $S_{x_i}^n=\{(x_1,x_2,...,x_{n+1}):x_1^2+x_2^2+...+x_{n+1}^2=1, x_i>0\}$ and $S_{-x_i}^n=\{(x_1,x_2,...,x_{n+1}):x_1^2+x_2^2+...+x_{n+1}^2=1, x_i<0\}$.
\end{center}  

Then we create a map $P_{x_i}:S_{x_i}^n \rightarrow D^n$ which projects these half-spheres into the unit disc. These maps are defined as\\ $P_{x_i}((x_1,x_2,...,x_{n+1}))=(x_1,x_2,...,\hat{x_i},...,x_{n+1})$. One can easily verify that each of these maps are diffeomorphisms so the proof will be omitted. When we consider all of these covers and their maps we notice that every neighborhood of $S^n$ is diffeomorphic to an open subset of $\mathbb{R}^n$. This means $S^n$ is a real manifold with a smooth structure.

\begin{definition}[Topological Abstract Manifolds]
    A subset M of a topological space X is said to be a \emph{smooth abstract manifold of dimension n} if M is second-countable, Hausdorff and each $x \in M$ has a neighborhood that is homeomorphic to an open subset of $\mathbb{R}^n$. 
\end{definition}

This definition instantly makes us wonder what smoothness means in the case of topological manifolds. This notion is explained in the following definition. 

\begin{definition}[Smoothness of Abstract Manifolds]
    A topological manifold of dimension m can be defined by $(M,\{ U_i \}_{i \in I}, \{ \phi_i \}_{i \in I})$ where M is subset of a topological space, $\{ U_i \}_{i \in I}$ is a family of covers of M, and $\{ \phi_i \}_{i \in I}$ is a family of homeomorphisms defined by $\phi_i:U_i \rightarrow \mathbb{R}^{n}$, this map takes $U_i$ onto an open subset of $\mathbb{R}^{n}$. We call a abstract manifold \emph{smooth} if for each of these coordinate maps we have that 
    \begin{center}
        $\phi_i\circ \phi_j^{-1}:\phi_j(U_i\cap U_j)\rightarrow \phi_i(U_i\cap U_j)$
    \end{center}
    is smooth.
\end{definition}

An example of a smooth abstract manifold will be introduced in the following section on Real Projective Spaces. 

\section{Real Projective Spaces}

In this section we introduce Real Projective Spaces, $\mathbb{RP^n}$,  and show that they are smooth manifolds. 

\quad Intuitively, Real Projective space of dimension $n$ can be thought of as the set of all lines through the origin of the Euclidean space of dimension $n+1$, i.e. 
\[\mathbb{RP}^n = \{\text{lines through the origin in } \mathbb{R}^{n+1}\}\]

\begin{figure}
    \centering
    \includegraphics[width = 6cm]{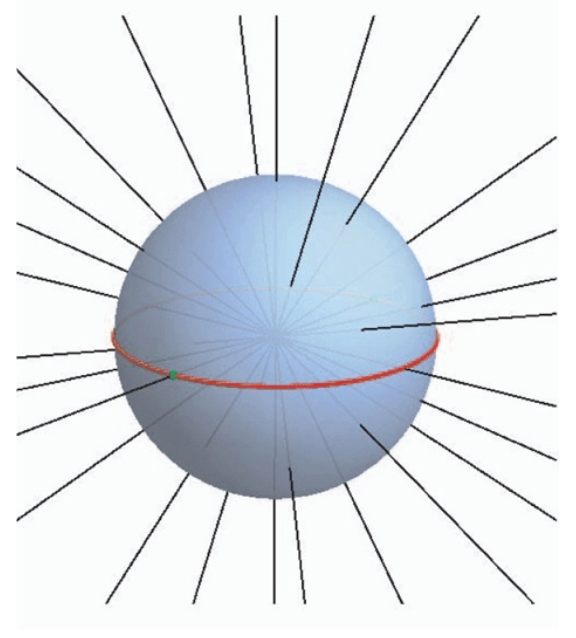}
    \caption{$\mathbb{RP}^2$ as a set of lines}
    \label{fig:enter-label}
\end{figure}

There are many rigorous identifications for $\mathbb{RP}^n$, but the one that will prove most useful for this topic is as follows:
\begin{definition}[Real projective spaces]
    \[\mathbb{RP}^n = \{[a_1,...,a_{n+1}] | (a_1...a_{n+1}) \neq 0 \in \mathbb{R}^{n+1}\}\]
    i.e. equivalence classes of vectors in $\mathbb{R}^{n+1}$, with the zero vector excluded. These equivalence classes are defined by relation of scalar multiplication.
\end{definition}

\newcommand{\bigslant}[2]{{\raisebox{.2em}{$#1$}\left/\raisebox{-.2em}{$#2$}\right.}}
We now show that $\mathbb{RP}^n$ is a smooth, abstract manifold.
First we must show that $\mathbb{RP}^n$ is a topological space. One can see that this is the case since 
\begin{center}
    $\mathbb{RP}^n=\bigslant{S^n}{(x\sim-x)}$
\end{center}
 and we have established that $S^n$ is a topological space so we know $\mathbb{RP}^n$ is endowed with the \emph{quotient topology}.

To see that $\mathbb{RP}^n$ is a manifold we define n+1 open covers in a similar way to the $S^n$ example. First we let $U_i,U_j \subset \mathbb{RP}^n$ with $U_i, U_j$ open and $i<j$. These covers are defined as:

\begin{itemize}
    \item[] $U_i=\{[x_1,x_2,...,x_i,...,x_{n+1}]|x_i\neq0\}=\{[\frac{x_1}{x_i},\frac{x_2}{x_i},...,1,...,\frac{x_{n+1}}{x_i}]|x_i \neq 0\}$
    \item[] $U_j=\{[x_1,x_2,...,x_j,...,x_{n+1}]|x_j\neq0\}
    =\{[\frac{x_1}{x_j},\frac{x_2}{x_j},...,1,...,\frac{x_{n+1}}{x_j}]|x_j \neq 0\}$
\end{itemize}
Now we define the coordinate maps of these covers, along with their inverses.

\begin{itemize}
    \item[] $\phi_i:U_i\rightarrow \mathbb{R}^n$, defined by $[x_1,x_2,...,x_i,...,x_{n+1}] \mapsto (\frac{x_1}{x_i},\frac{x_2}{x_i},...,\hat{\frac{x_i}{x_i}},...,\frac{x_{n+1}}{x_i})$ 
    \item[] $\phi_j:U_j\rightarrow \mathbb{R}^n$, defined by $[x_1,x_2,...,x_j,...,x_{n+1}] \mapsto (\frac{x_1}{x_j},\frac{x_2}{x_j},...,\hat{\frac{x_j}{x_j}},...,\frac{x_{n+1}}{x_j})$
\end{itemize}

the inverses are: 

\begin{itemize}
    \item[] $\phi_i^{-1}:\mathbb{R}^n \rightarrow U_i$, defined by $(a_1,a_2,...,a_n) \mapsto [a_1,a_2,...a_{i-1},1,a_{i},...,a_n]$
    \item[] $\phi_j^{-1}:\mathbb{R}^n \rightarrow U_j$, defined by $(a_1,a_2,...,a_n) \mapsto [a_1,a_2,...a_{j-1},1,a_{j},...,a_n]$
\end{itemize}
Now we compose the maps $\phi_i\circ\phi_j^{-1}$ and $\phi_j\circ\phi_i^{-1}$, which result in the following:
\begin{itemize} 
    \item[] $\phi_i \circ \phi_j^{-1}((a_1,a_2,...,a_n))=(\frac{a_1}{a_i},\frac{a_2}{a_i},...,\hat{\frac{a_i}{a_i}},...,\frac{a_{j-1}}{a_i},\frac{1}{a_i},\frac{a_j}{a_i},\frac{a_{j+1}}{a_i},...,\frac{a_n}{a_i})$
     \item[] $\phi_j \circ \phi_i^{-1}((a_1,a_2,...,a_n))=(\frac{a_1}{a_j},\frac{a_2}{a_j},...,\frac{a_{i-1}}{a_j},\frac{1}{a_j},\frac{a_i}{a_j},\frac{a_{i+1}}{a_j},..., \frac{a_{j-1}}{a_j},\hat{\frac{a_j}{a_j}},\frac{a_{j+1}}{a_j},...,\frac{a_n}{a_j})$
\end{itemize}
All of $\mathbb{RP}^n$ is covered and each of these covers has  a coordinate map associated with it. These covers and their maps satisfy the smoothness condition stated in definition 2.10. This means that $\mathbb{RP}^n$ is a smooth abstract manifold 

\section{Schwartz's Observation}
Schwartz functions are functions $f: \mathbb{R}^n \rightarrow \mathbb{R}$ such that $f$ is smooth and has 'rapid decay' at $\infty$. We formalize this definition as
\begin{definition}[Schwartz Property]
    A function $f: \mathbb{R}^n \rightarrow \mathbb{R}$ is Schwartz if:
    \begin{itemize}
        \item[] $f$ is smooth
        \item[] $\forall \alpha, \beta \in \mathbb{Z}_{\geq 0}^n, \quad \sup|x^{\alpha} \frac{\partial^{\beta}}{\partial x^{\beta}} f(x)| < \infty$
    \end{itemize}
\end{definition}

Schwartz observed that a Schwartz function on $\mathbb{R}^n$ extends to a smooth function on $S^n$ and is flat, or vanishes, at the infinity point. The inverse is also true: if $g$ is a smooth function on $S^n$ that is flat at $\{\infty\}$, then the restriction of $g$ to the reals, $g|_{\mathbb{R}^n}$, is a Schwartz function. \\

One classic example of a Schwartz function is $f(x) = e^{-x^2}$. 
\begin{proof}{$e^{-x^2}$ is Schwartz}
\begin{itemize}
    \item[] First we notice that $e^{-x^2}$ is symmetric about $y=0$.
    \item[] Now we show that $\lim_{x\to\infty}x^{n}e^{-x^2}$ goes to zero for all $n\in \mathbb{N}$. We rewrite $\lim_{x\to\infty}x^{n}e^{-x^2}$ as $\lim_{x\to\infty} \frac{x^{n}}{e^{x^2}}$ and apply L'Hopital's rule. This gives us $\lim_{x\to\infty} \frac{nx^{n-1}}{2xe^{x^2}}=\lim_{x\to\infty} \frac{n}{2}\frac{x^{n-2}}{e^{x^2}}$. Repeating this process $\frac{n}{2}$ times we get $\lim_{x\to\infty} \frac{(n)(n-2)...(n-2(\frac{n}{2}))}{2^{\frac{n}{2}}}\frac{1}{e^{x^2}}$. Clearly this goes to zero.
    \item[] Since differentiating $e^{-x^2}$ results in $e^{-x^2}$ multiplied by $-2x$ we know due due the power rule that every derivative of $e^{-x^2}$ is of the form $P(x)e^{-x^2}$ where $P(x)$ is a polynomial. 
    \item[] Thus, $e^{-x^2}$ is Schwartz.
\end{itemize}
\end{proof}
The higher dimensional version of this equation, $g(x_1, ..., x_n) = e^{-(x_1^2+...+x_n^2)}$, is still Schwartz. In fact, all Gaussian functions are Schwartz. An interesting example of a function that is \emph{not} Schwartz is given below.
\begin{proof}{$h(x) = e^{-x^2}\sin{e^{x^2}}$ is \textit{not} Schwartz}
    \begin{itemize}
        \item[] We notice this function is still symmetric about the y axis, but it exhibits some oscillatory behavior before settling to 0 as $x \rightarrow \infty$.
        \item[] Of course, we also need $h'(x)$ and all other derivatives to vanish at infinity. Our problem arises at the first derivative, \\
        $h'(x) = -2xe^{-x^2}\sin{e^{x^2}} + 2xe^{2x^2}\cos{e^{x^2}}$. 
        \item[] $\lim_{x\to\infty}-2xe^{-x^2}\sin{e^{x^2}} + 2xe^{x^4}\cos{e^{x^2}}$ is clearly not bounded, as the second term blows up to infinity as x goes to infinity. 
        \item[] Thus, $e^{-x^2}\sin{e^{x^2}}$ is not Schwartz. 
    \end{itemize}
\end{proof}

This concept of Schwartz functions and their behavior at the 'infinity point' gives us the following lemma, which will be used to show our end result:

\begin{lemma}[]
     Let $g(t): \mathbb{R}^n\text{\textbackslash}\{t_1=0\} \rightarrow \mathbb{R}$ be a smooth function. Then
        $g$ extends to a "flat" function at $t_1=0 \iff$ $\forall a \in \{t_1=0\}$, $\exists r>0$, $\forall\alpha\in \mathbb{Z}^n_{\geq0}$, $\forall p \in \mathbb{Z}_{\geq 0}$,\\ 
        \[\sup\limits_{t\in B_r(a)\text{\textbackslash}\{t_1=0\}} |t_1^{-p} \frac{d^{\alpha} g}{dt^{\alpha}}| < \infty\]
\end{lemma}

In this lemma, \textit{flatness} means that the function and all its derivatives go to zero at infinity.

We now aim to extend this to establish 
$\mathbb{RP}^n = \mathbb{R}^n \sqcup \mathbb{RP}^{n-1}$. 

\begin{figure}
    \centering
    \includegraphics[width = 8cm]{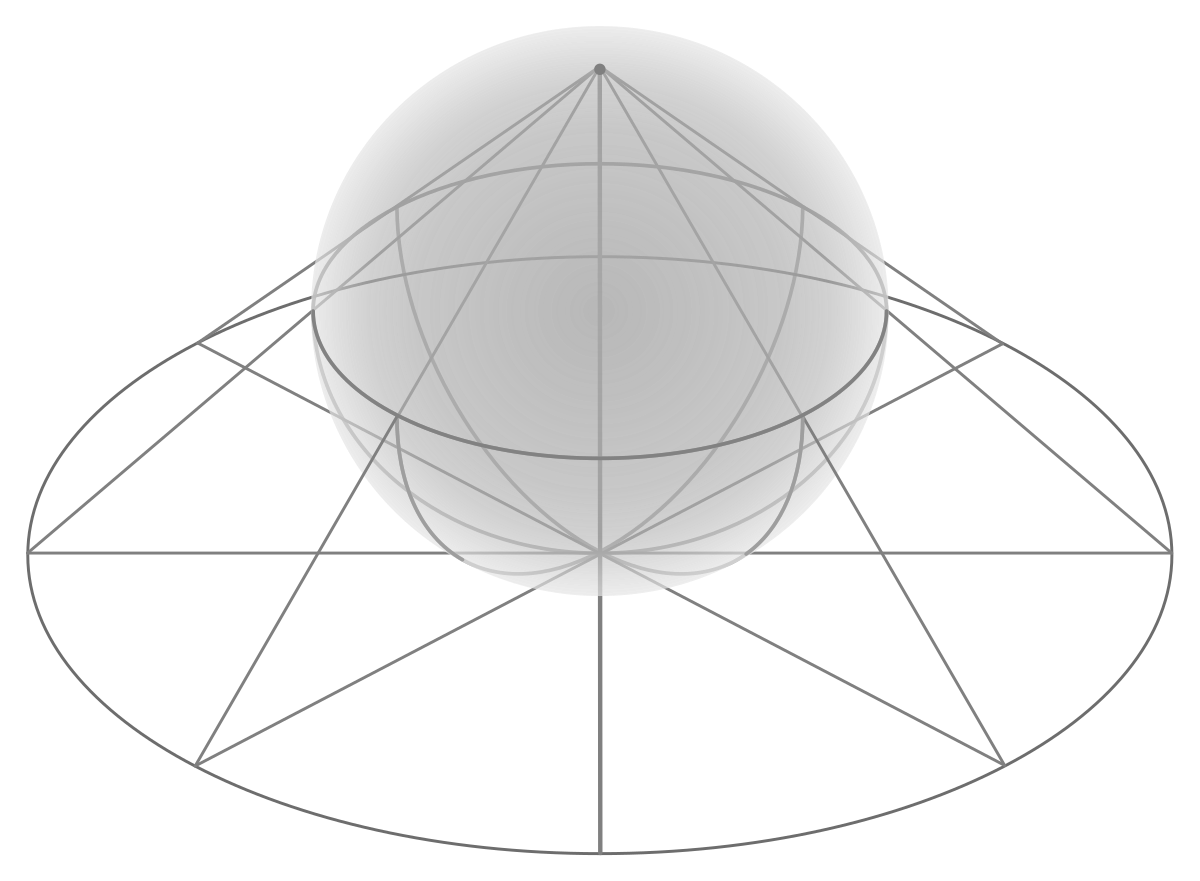}
    \caption{Stereographic Projection of $S^2$ onto $\mathbb{R}^2$}
    \label{fig:enter-label}
\end{figure}

\section{$\mathbb{RP}^n = \mathbb{R}^n \sqcup \mathbb{RP}^{n-1}$}
\quad We consider $\mathbb{RP}^n = U_i\sqcup V_i$ where $U_i \cong \mathbb{R}^n$ and $V_i \cong \mathbb{RP}^{n-1}$. We begin by showing that $U_i \cong \mathbb{R}^n$ and $V_i \cong \mathbb{RP}^{n-1}$. To do this we consider the open subset $U_i \subseteq \mathbb{RP}^n$ such that $U_i = \{[x_1, x_2,...,x_{n+1}] : x_i \neq 0\}$. Using the smooth coordinate map $\phi_i:U_i \rightarrow \mathbb{R}^n$ that takes $[x_1, ..., x_{n+1}] \mapsto (\frac{x_1}{x_i}, \frac{x_2}{x_i},..., \hat{\frac{x_i}{x_i}}, ... \frac{x_n}{x_i})$ we notice that $\phi_i(U_i)=\mathbb{R}^n$. 

If we consider $\mathbb{RP}^n$\textbackslash $U_i = \{x_1, x_2,...,x_i,...x_n : x_i = 0\}$, we see that we can identify this with $\mathbb{RP}^{n-1}$. Thus, we can define $\mathbb{RP}^n = U_i\sqcup \mathbb{RP}^{n-1}$
\\\\
Now, if we have a Schwartz function on $U_1 \subseteq \mathbb{RP}^n$, we want to show that it can extend to the remaining $V_i \cong \mathbb{RP}^{n-1}$. In order to get there, we need to consider the relations between $U_1, U_2$ on $\mathbb{RP}^n$ by exploring how the partial derivatives of their transition functions interact with each other. \\\\
We will begin with an example where n=3 to understand how the calculations work, and then we will extrapolate this to the n-dimensional case. \\
Consider a Schwartz function $f(s_1, s_2, s_3)$ on $\phi_1(U_1) \longleftrightarrow g(t_1, t_2, t_3)$ on $\phi_2(U_2)$. These functions are related by the following compositions:
\begin{align*}
    f\circ(\phi_1\circ\phi_2^{-1})(t_1, t_2, t_3) &= g(t_1, t_2,t_3)\\
    g\circ(\phi_2\circ\phi_1^{-1})(s_1, s_2, s_3) &= f(s_1, s_2, s_3)
\end{align*}
Where 
\begin{align*}
    \phi _1 \circ \phi _2 ^{-1}: \phi_2(U_1\cap U_2) & \to \phi _1 (U_1\cap U_2) \\
    (t_1, t_2, t_3) & \mapsto (s_1, s_2, s_3) = (\frac{1}{t_1}, \frac{t_2}{t_1}, \frac{t_3}{t_1})
\end{align*}
and 
\begin{align*}
    \phi _2 \circ \phi _1 ^{-1}: \phi_1(U_1\cap U_2) & \to \phi _2 (U_1\cap U_2)\\
    (s_1, s_2, s_3) & \mapsto (t_1, t_2, t_3) = (\frac{1}{s_1}, \frac{s_2}{s_1}, \frac{s_3}{s_1})
\end{align*}\\
Then we can consider the partial derivatives: 
\begin{align*}
    \frac{\partial f}{\partial s_1} &\longleftrightarrow \frac{\partial g}{\partial t_1}, \frac{\partial g}{\partial t_2}, \frac{\partial g}{\partial t_3}\\
    \frac{\partial g}{\partial t_1} &\longleftrightarrow \frac{\partial f}{\partial s_1}, \frac{\partial f}{\partial s_2}, \frac{\partial f}{\partial s_3}
\end{align*}

We would like to show that the partial derivatives of each transition function are also related to the derivatives of the inverse transition function. The calculations are as follows: 

\begin{itemize}
    \item[] $\frac{\partial g}{\partial t_1} = \frac{\partial (f(s_1, s_2, s_3))}{\partial t_1} = \frac{\partial s_1}{\partial t_1}\frac{\partial f}{\partial s_1} + \frac{\partial s_2}{\partial t_1}\frac{\partial f}{\partial s_2} + \frac{\partial s_3}{\partial t_1}\frac{\partial f}{\partial s_3} = -\frac{1}{t_1^2}\frac{\partial f}{\partial s_1} + t_2\frac{\partial f}{\partial s_2} + t_3\frac{\partial f}{\partial s_3}$
    \item[] $\frac{\partial g}{\partial t_2} = \frac{\partial s_1}{\partial t_2}\frac{\partial f}{\partial s_1} + \frac{\partial s_2}{\partial t_2}\frac{\partial f}{\partial s_2} + \frac{\partial s_3}{\partial t_3}\frac{\partial f}{\partial s_3} = 0 + \frac{1}{t_1}\frac{\partial f}{\partial s_2}+0 = \frac{1}{t_1}\frac{\partial f}{\partial s_2}$
    \item[] $\frac{\partial g}{\partial t_3}= \frac{\partial s_1}{\partial t_3}\frac{\partial f}{\partial s_1} + \frac{\partial s_2}{\partial t_3}\frac{\partial f}{\partial s_2} + \frac{\partial s_3}{\partial t_3}\frac{\partial f}{\partial s_3} = 0 + 0 + \frac{1}{t_1}\frac{\partial f}{\partial s_3} = \frac{1}{t_1}\frac{\partial f}{\partial s_3}$
\end{itemize}

If we solve each of the equivalences for their respective $\frac{\partial f}{\partial s_n}$, we get:

\begin{itemize}
    \item[] $\frac{\partial f}{\partial s_2} = t_1\frac{\partial g}{\partial t_3}$
    \item[] $\frac{\partial f}{\partial s_2} = t_1\frac{\partial g}{\partial t_2}$
    \item[] $\frac{\partial f}{\partial s_1} = -t_1^2\frac{\partial g}{\partial t_1}-t_1t_2\frac{\partial g }{\partial t_2}- t_1t_3\frac{\partial g}{\partial t_3}$    
\end{itemize}

Finally, arranging these solutions in matrix form, we obtain 

\[ \begin{bmatrix}
    \partial s_1\\
    \partial s_2\\
    \partial s_3
    \end{bmatrix} = 
    \begin{bmatrix}
    -t_1^2 & -t_1t_2 & -t_1t_3\\
    0 & t_1 & 0 \\
    0 & 0 & t_1\\
    \end{bmatrix}
    \begin{bmatrix}
        \partial t_1\\
        \partial t_2\\
        \partial t_3
    \end{bmatrix}\]
We can follow the same process to obtain relations in terms of $\frac{\partial g}{\partial t_n}$, resulting in the following matrix: 
\[ \begin{bmatrix}
    \partial t_1\\
    \partial t_2\\
    \partial t_3
    \end{bmatrix} = 
    \begin{bmatrix}
    -s_1^2 & -s_1s_2 & -s_1s_3\\
    0 & s_1 & 0 \\
    0 & 0 & s_1\\
    \end{bmatrix}
    \begin{bmatrix}
        \partial s_1\\
        \partial s_2\\
        \partial s_3
    \end{bmatrix}\]

\quad Each of our covers is related to the other by these transformations described in the matrices. Thus, we can obtain $t_1^{-\beta} \frac{d^{\alpha} g}{dt^{\alpha}}$ as a sum of polynomials in terms of $s$ and vice versa as $s_1^{-\beta} \frac{d^{\alpha} f}{ds^{\alpha}}$ in terms of $t$. This will result in equivalences of the following form: 

\[t_1^{-\beta_1}\frac{\partial^\alpha g}{\partial t^{\alpha}}(t) = \sum_{|\beta| < |\alpha|} s_1^{\beta_1} P_{\beta}(s) \frac{\partial^\beta f}{\partial s^{\beta}}(s)\]
\[s_1^{-\beta_1}\frac{\partial^\alpha f}{\partial s^{\alpha}}(s) = \sum_{|\gamma| < |\alpha|} t_1^{\beta_1} P_{\gamma}(t) \frac{\partial^\gamma g}{\partial t^{\gamma}}(t) \]

\noindent where $P_{\beta}(s)$ is some polynomial in terms of $s$. Because we know each $\frac{\partial^\beta f}{\partial s^{\beta}}(s)$ go to 0 as Schwartz functions, and these polynomials are bounded, then we can conclude that $\sum_{|\beta| < |\alpha|} s_1^{\beta_1} P_{\beta}(s) \frac{\partial^\beta f}{\partial s^{\beta}}(s)$ is bounded, and thus $t_1^{-\beta_1}\frac{\partial^\alpha g}{\partial t^{\alpha}}(t)<\infty$. Thus by Lemma 3.1, we have that $g$ extends to a flat function at $t_1 = 0$. \qedsymbol\\

This observation holds for $\mathbb{RP}^n$. We will proceed using the same covers and transition function defined in section 2.3 and techniques used for the $\mathbb{RP}^3$ case. We let $i<j$ and take our open covers to be
\begin{itemize}
    \item[] $U_i = \{[x_1, x_2,...,x_{n+1}] : x_i \neq 0\}$  
    \item[] $U_j =\{[x_1, x_2,...,x_{n+1}] : x_j \neq 0\}$
\end{itemize}
with smooth functions and their inverses:
\begin{itemize}
    \item[] $\phi_i:U_i\rightarrow \mathbb{R}^n$, defined by $[x_1,x_2,...,x_i,...,x_{n+1}] \mapsto (\frac{x_1}{x_i},\frac{x_2}{x_i},...,\hat{\frac{x_i}{x_i}},...,\frac{x_{n+1}}{x_i})$ 
    \item[] $\phi_j:U_j\rightarrow \mathbb{R}^n$, defined by $[x_1,x_2,...,x_j,...,x_{n+1}] \mapsto (\frac{x_1}{x_j},\frac{x_2}{x_j},...,\hat{\frac{x_j}{x_j}},...,\frac{x_{n+1}}{x_j})$
    \item[] $\phi_i^{-1}:\mathbb{R}^n \rightarrow U_i$, defined by $(a_1,a_2,...,a_n) \mapsto [a_1,a_2,...a_{i-1},1,a_{i},...,a_n]$
    \item[] $\phi_j^{-1}:\mathbb{R}^n \rightarrow U_j$, defined by $(a_1,a_2,...,a_n) \mapsto [a_1,a_2,...a_{j-1},1,a_{j},...,a_n]$
\end{itemize}

The transition functions are given by:
\begin{itemize} 
    \item[] $\phi_i \circ \phi_j^{-1}((a_1,a_2,...,a_n))=(\frac{a_1}{a_i},\frac{a_2}{a_i},...,\hat{\frac{a_i}{a_i}},...,\frac{a_{j-1}}{a_i},\frac{1}{a_i},\frac{a_j}{a_i},\frac{a_{j+1}}{a_i},...,\frac{a_n}{a_i})$
     \item[] $\phi_j \circ \phi_i^{-1}((a_1,a_2,...,a_n))=(\frac{a_1}{a_j},\frac{a_2}{a_j},...,\frac{a_{i-1}}{a_j},\frac{1}{a_j},\frac{a_i}{a_j},\frac{a_{i+1}}{a_j},..., \frac{a_{j-1}}{a_j},\hat{\frac{a_j}{a_j}},\frac{a_{j+1}}{a_j},...,\frac{a_n}{a_j})$
\end{itemize}

Now, just as before, we must find the relationship between the derivatives of these transition functions. This results in the following matrices: \\
In terms of $\frac{\partial f}{\partial s_n}$
\[ \begin{bmatrix}
    \partial_{s_1} \\
    \partial_{s_2}\\
    \partial_{s_3}\\
    ...\\
    \partial_{s_n}
    \end{bmatrix} = 
    \begin{bmatrix}
    -t_1^2 & -t_1t_2 & -t_1t_3 & ... & -t_1t_n\\
    0 & t_1 & 0 &...&0\\
    0 & 0 & t_1&...&0\\
    ...&...&...&...&...\\
    0&0&0&...&t_1
    \end{bmatrix}
    \begin{bmatrix}
        \partial_{t_1}\\
        \partial_{t_2}\\
        \partial_{t_3}\\
        ...\\
        \partial_{t_n}
    \end{bmatrix}\]

and in terms of  $\frac{\partial g}{\partial t_n}$:
     \[ \begin{bmatrix}
    \partial_{t_1} \\
    \partial_{t_2}\\
    \partial_{t_3}\\
    ...\\
    \partial_{t_n}
    \end{bmatrix} = 
    \begin{bmatrix}
    -s_1^2 & -s_1s_2 & -s_1s_3 & ... & -s_1s_n\\
    0 & s_1 & 0 &...&0\\
    0 & 0 & s_1&...&0\\
    ...&...&...&...&...\\
    0&0&0&...&s_1
    \end{bmatrix}
    \begin{bmatrix}
        \partial_{s_1}\\
        \partial_{s_2}\\
        \partial_{s_3}\\
        ...\\
        \partial_{s_n}
    \end{bmatrix}\]

It is important to note that all previous calculations are done using the first two open subsets, $U_1$ and $U_2$. These calculations can be repeated with any transition functions between $U_1$ and $U_i$ for any $i>2$, which will still result in matrices for transition functions with polynomial entries, though they aren't nicely upper-triangular. Nonetheless, the Lemma can still be applied and our result will hold.

\section{Appendix}

Here, we include some additional theorems, definitions, and examples that are related to our topic, but did not warrant a spotlight. \\

\begin{theorem}[Heine-Borel]
Let $S \subset \mathbb{R}^n$. Then
\begin{center}
    $S$ is closed and bounded $\iff S$ is compact
\end{center}
\end{theorem}

As our work is focused on finding a compactification of $\mathbb{RP}^n$, which is a very specific case. We add the more general definition for this concept below:
\begin{definition}[Compactification]
    Compactification is the mathematical process of embedding topological spaces into compact Hausdorff spaces.  
\end{definition}
\quad One exercise explored by the authors was the stereographic projection of $S^n$ onto $R^n$. We derived the equation for stereographic projection as well as its inverse below:\\
Given a point $(y, x_1, ...,x_n)\in S^n$ and the north pole of $S^n$ $(1,0,...,0)$ we find where each point lands when y goes to zero. we find this point by solving $0=\frac{y-1}{x_i}x+1$ for x for each $x_i$. We get that $x=\frac{x_i}{1-y}$. This tells us that $f((y, x_1,...,x_n))=(\frac{x_1}{1-y},\frac{x_2}{1-y},...,\frac{x_n}{1-y})$.\\
\\
We notice that the distance of $f(\vec{x})$ from $\vec{0}$ is completely determined by the y component of $\vec{x}$. So we let $d=(\frac{x_1}{1-y})^2+(\frac{x_2}{1-y})^2+...(\frac{x_n}{1-y})^2=\frac{x_1^2+x_2^2+...+x_n^2}{(1-y)^2}$ but we know $x_1^2+x_2^2+...+x_n^2=1-y^2$ since $(x_1,y,x_2,...,x_n)$ lies on $S^n$. This means $d=\frac{1-y^2}{(1-y)^2}=\frac{1+y}{1-y}$. From here we can find the value of y. Doing some calculations we find that $y=\frac{d-1}{d+1}$ This means $f^{-1}:\mathbb{R}^n \rightarrow S^n$ is defined by $f^{-1}((x_1,x_2,...,x_n))=((1-\frac{d-1}{d+1})x_1,\frac{d-1}{d+1},(1-\frac{d-1}{d+1})x_2,...,(1-\frac{d-1}{d+1})x_n)=(\frac{2x_1}{d+1},\frac{d-1}{d+1},\frac{2x_2}{d+1},...,\frac{2x_n}{d+1})$.\\

\section{Acknowledgements}
We would like to thank Dr. Ahmad Reza Haj Saeedi Sadegh for leading this project, as well as his guidance and expertise throughout this project. I would like to thank J.D for motivating me to write something which I could dedicate to them. A special thanks to A.D.S. who always inspired diligence. We would also like to thank the Northeastern University College of Science, the Northeastern University Department of Mathematics, and the NSF-RTG grant "Algebraic Geometry and Representation Theory at Northeastern University" (DMS-1645877).

\nocite{*}
\printbibliography

\end{document}